\documentclass[12pt]{amsart}
\usepackage[T1]{fontenc}
\usepackage{mathptmx}

\usepackage{comment}
\usepackage[colorlinks=true,linkcolor=magenta,citecolor=blue]{hyperref}     
\numberwithin{equation}{section}
\usepackage{mathtools}
\usepackage{amsmath, amssymb,amsthm}
\usepackage[capitalise]{cleveref}   
\usepackage{array}
\usepackage{graphicx}
\usepackage{float}
\usepackage{caption,subcaption}
\usepackage{tikz}
\usepackage{tikz-cd}
\usetikzlibrary{matrix,arrows,arrows.meta}
\usepackage[mathscr]{euscript}
\usepackage{xcolor}

\def\Abb{\mathbb{A}} 
\def\Cbb{\mathbb{C}}

\def\Kbb{\mathbb{K}} 
 
 \def\Pbb{\mathbb{P}}
\def\Qbb{\mathbb{Q}} 
\def\Sbb{\mathbb{S}}

 \def\Zbb{\mathbb{Z}}

\def\CCal{\mathcal{C}}  
\def\Ecal{\mathcal{E}} 
 \def\Hcal{\mathcal{H}}
 
 \def\Lcal{\mathcal{L}}
 
\def\Ocal{\mathcal{O}}

 \def\Vcal{\mathcal{V}}

\newtheorem{thm}{Theorem}[section]
\newtheorem{mainthm}{Theorem}

\newtheorem{lem}[thm]{Lemma}
\newtheorem{ex}[thm]{Example}
\theoremstyle{definition}
\newtheorem{rem}[thm]{Remark}

\theoremstyle{plain}

\newtheorem{prop}[thm]{Proposition}

\crefname{thm}{Theorem}{Theorems}
\Crefname{thm}{Theorem}{Theorems}
\crefname{lem}{Lemma}{Lemmas}
\Crefname{lem}{Lemma}{Lemmas}
\crefname{prop}{Proposition}{Propositions}
\Crefname{prop}{Proposition}{Propositions}
\crefname{cor}{Corollary}{Corollaries}
\Crefname{cor}{Corollary}{Corollaries}
\crefname{defi}{Definition}{Definitions}
\Crefname{defi}{Definition}{Definitions}
\crefname{rem}{Remark}{Remarks}
\Crefname{rem}{Remark}{Remarks}

\newcommand\an{\mathrm{an}}

\newcommand\Jac{\mathrm{Jac}}

\newcommand\gal{\operatorname{Gal}}

\DeclareMathOperator{\ddc}{dd^c}

\DeclareMathOperator{\id}{id}

\DeclareMathOperator{\Fix}{Fix}

\DeclareMathOperator{\Spec}{Spec}

\DeclareMathOperator{\Aut}{Aut}

\begin{document}
\title{The Dynamical Manin-Mumford Conjecture for automorphisms of H\'enon type}

\author{Thomas Gauthier}
\address{Laboratoire de Math\'ematiques d'Orsay, B\^atiment 307, Universit\'e Paris-Saclay, 91405 Orsay Cedex, France}
\email{thomas.gauthier@universite-paris-saclay.fr}
\thanks{The first author's research is partially supported by the Institut Universitaire de France. The second author's research is partially supported by the programme ATRACT of Région Bourgogne-Franche-Comté under the project ADYAUS (RECH-ATRAC-000012) and the French National Research Agency under the project DynAtrois (ANR-24-CE40-1163)}
\author{Yugang Zhang}
\address{Yugang Zhang, Université Bourgogne Europe, CNRS, IMB UMR 5584, 21000 Dijon, France}
\email{yugang.zhang@ube.fr}

\begin{abstract}
    We prove the dynamical Manin--Mumford conjecture for polynomial automorphisms of the complex affine plane of dynamical degree $>1$ proposed by Dujardin and Favre. More precisely, we prove that an irreducible algebraic curve contains infinitely many periodic points of such an automorphism if and only if it is the fixed locus of a polynomial involution conjugating some iterate of the automorphism to its inverse.
\end{abstract}
\maketitle

\section{Introduction}

The classical Manin--Mumford theorem, proved by Raynaud~\cite{Raynaud83,Raynaud83Courbes}, characterizes the irreducible subvarieties of a complex abelian variety containing a Zariski dense set of torsion points: they are precisely the translates of abelian subvarieties by torsion points. 
Since torsion points of an abelian variety are precisely the preperiodic points of the multiplication map $[n]$, for all $n\geq2$, Raynaud's theorem can be interpreted as a dynamical statement. Recall that, given a self-map $f$ of a set $X$, a point $x\in X$ is \emph{preperiodic} if its forward orbit $\{f^n(x)\}_{n\geq0}$ is finite, and \emph{periodic} if some positive iterate fixes it. When $f$ is a bijection, these notions coincide.

This analogy led Zhang~\cite{Zhang95} to propose a dynamical Manin--Mumford conjecture~\cite{Zhang95}: a subvariety containing a Zariski dense set of preperiodic points of a polarized endomorphism should itself be preperiodic. Counterexamples of Ghioca--Tucker--Zhang \cite{GTZ11} show that this formulation fails in general. Modified versions were subsequently proposed in \cite{GTZ11,GT21}, while DeMarco--Mavraki \cite{DM24} formulated a relative version within families. Roughly speaking, one expects that the original conjecture holds unless the dynamical system and the subvariety are very special. Although the conjecture remains widely open, substantial progress has been made for split endomorphisms~\cite{MS14,GNY1,GNY2,MSW23}, lifts of Frobenius~\cite{Xie18}, and regular polynomial endomorphisms of $\Cbb^2$~\cite{DMM-DFR}.

\medskip

In this text we consider the dynamical Manin--Mumford problem for polynomial automorphism of $\Abb^2_\Cbb$. For such $f$, its \emph{dynamical degree} is $\lambda_f:=\lim \deg(f^n)^{1/n} $, where $\deg$ denotes the algebraic degree, i.e. the maximal degrees of the components of $f$. 
By the classification of Friedland--Milnor~\cite{FM89}, if $\lambda_f=1$, then the set of periodic points is algebraic, hence the dynamical Manin--Mumford problem is not interesting. 

A polynomial automorphism of $f$ with $\lambda_f>1$ is called \emph{of H\'enon type}. Such an automorphism is \emph{reversible} if it is conjugate to its inverse by an involution $\sigma$ possessing a curve $C$ of fixed points; see G\'omez--Meiss~\cite{GM03,GM04}. For instance, $f(x,y)=(y,y^2-x)$ is reversible with respect to the involution $\sigma(x,y)=(y,x)$, whose fixed locus is the diagonal $C=\{x=y\}$. In this case there are infinitely many periodic points on $C$~\cite[Proposition 7.1]{DF-MM}.
Dujardin--Favre conjectured that every such curve arises in this way, up to iterates of $f$. As a first step, they proved that if there exist infinitely many periodic points on an algebraic curve $C$, then the Jacobian $\mathrm{Jac}(f)$ is algebraic and all its Galois conjugates have complex modulus one~\cite[Theorem A]{DF-MM}. 

Our main theorem confirm their conjecture.

\begin{mainthm}\label{thm:main}
    Let $f$ be a polynomial automorphism of $\Cbb^2$ of H\'enon type. Let $C\subset \Cbb^2$ be an irreducible algebraic curve. 
    Then $C$ contains infinitely many $f$-periodic points if and only if there exists an integer $r\ge 1$ and a polynomial involution $\sigma$ of $\Cbb^2$ such that its fixed points are $C$ and $\sigma \circ f^{r} \circ \sigma ^{-1}=f^{-r}$.
\end{mainthm}
Note that over number fields we prove a stronger Bogomolov-type result, see Theorem~\ref{mainthm:nbfield} below.

\medskip

The strategy of the proof decomposes in three main steps. For simplicity assume that the Zariski closure of $C$ in $\Pbb^2$ is smooth. 
First, we let $F=f\times f^{-1}:\Abb^2\times\Abb^2\to\Abb^2\times\Abb^2$ and $\Delta_C=\{(c,c):c\in C\}$. Using the fact that $f$ admits no periodic algebraic curve~\cite{BS91}, the first step of the proof reduces the theorem to showing that
\[
    \bigcup_{n\in\Zbb}F^n(\Delta_C)
       \quad\text{is not Zariski dense in }\Abb^2\times\Abb^2.
\]
This is done in Section~\ref{sect:reduction}.

\medskip

In the second step, we first assume that $f$ and $C$ are defined over a number field $\mathbb{K}$. Using the arithmetic equidistribution of periodic points of Yuan~\cite{Yuan} or Thuillier~\cite{Thuillier} and the arithmetic Hodge index theorem of Yuan--Zhang~\cite{Hodge-YZ}, together with a rigidity result of Cantat--Dujardin~\cite{CD24}, we show that
\[
G^+_{f,v}=\alpha G^-_{f,v}
\]
on $C_v^\mathrm{an}$ for any place $v\in M_\mathbb{K}$, where $\alpha>0$ is a rational constant independent of the place; see Proposition~\ref{prop-equidistrib}.

The proportionality of Green functions implies an analytic relation at each branch at infinity of the curve via B\"ottcher coordinates. More precisely, Proposition~\ref{prop-bottcher} gives integers $A,B\ge 1$, and another sufficiently large integer $N_0\ge 1$ such that
\begin{align}
    \left(\phi_f^+\circ f^{N_0}\right)^A
       =\left(\phi_f^-\circ f^{-N_0}\right)^B
\end{align}
on a punctured neighborhood of a branch at infinity of $C$.  This is done in Section~\ref{sect:bottcher}.
In Section~\ref{sect:bottcher-complex}, we extend this B\"ottcher relation to those $f$ and $C$ defined over $\Cbb$ relying on a specialization argument; see Proposition~\ref{prop:complex-bottcher}.

Finally, for the last step, we use Xie's algebraization theorem~\cite{Xie2024} to deduce the required
non-Zariski density of $\bigcup_{n\in\Zbb}F^n(\Delta_C)$ from this B\"ottcher relation. We do this over the function field of the curve $C$. This completes the proof. This is done in Section~\ref{sect:xie}.

For comparison, Dujardin--Favre use arithmetic equidistribution to obtain that $G^+_{f,v}=\alpha G^-_{f,v}+H_v$ on  $C_v^{\mathrm{an}}$, where $\alpha>0$ and $H_v$ are harmonic on $C_v^\mathrm{an}$. At Archimedean places, they then use Pesin theory and laminar currents
to compute the local lower H\"older exponents of the Green functions restricted to $C$ in terms of the Lyapunov exponents of the equilibrium
measure, and deduce that $|\Jac(f)|_v=1$.

Note also that this is not the first use of an algebraization result in the study of unlikely intersections problems in arithmetic dynamics: this kind of strategy was already used in \cite{FG,Xie2024}.

\subsection*{AI statement} No LLM was used at any stage of this work, in any way whatsoever.

\section{Reduction to a non-Zariski density property}\label{sect:reduction}
Throughout this section, let $k$ be a field of characteristic zero, and let $f\in\Aut(\Abb^2_k)$ be of H\'enon type. Denote by $\bar k$ an algebraic closure of $k$. We first recall the following geometric properties, due to Bedford, Lyubich and Smillie.

\begin{prop}\label{prop:basic-bs}
\begin{enumerate}
    \item There is no $f$-periodic irreducible algebraic curve in $\Abb^2_{\bar k}$; see \cite[Proposition~4.2]{BS91}.
    \item The periodic points of $f$ in $\Abb^2(\bar k)$ are Zariski dense in $\Abb^2_{\bar k}$; see \cite[Theorem~1]{BLS93}.
\end{enumerate}
\end{prop}

Set $X=\Abb^2_k$. Let $F:=f\times f^{-1}$ be the product automorphism
of $X\times X$. Let $Y$ be a geometrically irreducible proper closed subvariety of $X\times X$. Suppose that $Y$ is $F$-periodic, i.e. $F^r(Y)=Y$ for some $r>0$. Let $\pi_i$, $i=1,2$, be the restrictions to $Y$ of the two projections $X\times X\to X$.
    
\begin{lem}\label{lem:dim_2}
    Suppose $k$ is algebraically closed. The dimension of $Y$ is either zero or two.
\end{lem}
\begin{proof}
    Suppose first that
    $\dim Y=3$. Then $\pi_i(Y)$ is Zariski dense in $X$. Indeed, otherwise
    $\overline{\pi_i(Y)}$ is a curve or a point. 
    By Proposition~\ref{prop:basic-bs} (1), it can not be a curve. If it was a point, then
    $\dim Y\leq 2\neq3$. In any case we have a contradiction.

    There exists a nonempty open subset $U\subset X$ such that the fiber $Y_x:=\{y\in X \mid (x,y)\in Y\}$ is a curve for every $x\in U$. By Proposition~\ref{prop:basic-bs} (2), there exists $x\in U$ such that
    $f^r(x)=x$ up to replacing $r$ by a multiple. Thus $f^{r}(Y_x)=Y_x$, which is again contradiction. Therefore $\dim Y\le 2$.

    Suppose now that $\dim Y=1$. Again since $f$ admits no periodic curve, both projections are constant, which forces $Y$ to be a point, a contradiction.
\end{proof}

\begin{lem}\label{lem:graph}
    Suppose that $k$ is algebraically closed and both projections $\pi_i$ are dominant. Then there exist an integer $r\ge 1$ and a polynomial automorphism $\sigma$ of $X$ such that
    \[
        Y=\operatorname{Graph}(\sigma),
        \qquad
        \sigma \circ f^r=f^{-r}\circ \sigma.
    \]
\end{lem}
\begin{proof}
Put $h_1=f^r$ and $h_2=f^{-r}$. The projection $\pi_i$ is quasi-finite. Indeed, the set of points, over which its fibers have positive dimension, is finite and $h_i$-invariant. Thus we would have a $h_i$-periodic curve if nonempty.

Let $\widetilde{Y}$ be the normalization of $Y$, and let $g$ be the automorphism induced by $F^r$.
The induced projections $p_i:\widetilde{Y}\to X$ satisfy $p_i\circ g=h_i\circ p_i.$

It suffices to show that $p_i$ is an isomorphism. Indeed, then we can set $\sigma:=p_2\circ p_1^{-1}$. Since the image of $(p_1,p_2)$ is $Y$, we have $Y=\operatorname{Graph}(\sigma)$ and $\sigma\circ  f^r=f^{-r}\circ \sigma.$

Now we show that $p_i$ is an isomorphism. By Zariski's Main Theorem \cite[Tag 05K0]{stacks-project},
each $p_i$ factors as
\[
\widetilde{Y}\hookrightarrow Z_i
\overset{q_i}{\longrightarrow}X,
\]
where $q_i$ is finite and $Z_i$ is normal and irreducible. 

Put $E_i=Z_i\setminus\widetilde Y$. Since $\widetilde Y$ is affine, $E_i$ is either empty or has pure codimension one.
Likewise, purity of the branch locus \cite[Tag 0BMB]{stacks-project} shows that the branch locus
of $q_i$ is either empty or a union of curves. Both $q_i(E_i)$ and this branch locus are $h_i$-invariant.
Proposition~\ref{prop:basic-bs} shows that both are empty. Hence $p_i=q_i$ is finite and \'etale, thus isomorphism.
\end{proof}

Let $C$ be a geometrically irreducible algebraic curve of $\Abb^2_k$ and put $C_{\bar k}:=C\times_k \bar k$ the base change.
Denote by $\Delta_{C_{\bar k}}:= \{(c,c)\in\Abb_{\bar k}^2\times \Abb_{\bar k}^2 \mid c\in C \}$ the diagonal of $C_{\bar k}$.
\begin{lem}\label{lem:keylemma}
Suppose $\bigcup_{n\in\Zbb} F^n(\Delta_{C_{\bar k}})$ is not Zariski dense in $\Abb^2_{\bar k}\times\Abb^2_{\bar k}$. Then there exist $r\geq1$ and a polynomial involution $\sigma$ of $\Abb^2_k$ defined over $k$ such that
\[
    \Fix(\sigma)=C,
    \qquad
    \sigma\circ f^r=f^{r}\circ\sigma.
\]
\end{lem}

\begin{proof}
The Zariski closure of $\bigcup_{n\in\Zbb} F^n(\Delta_{C_{\bar k}})$ is fixed by $F$.
Let $Y$ be an irreducible component containing $\Delta_{C_{\bar k}}$. Then $f^r(Y)=Y$ for some $r\geq1$. Both projections of $Y$ are dominant: otherwise, the Zariski closure of one of their images would give an $f$-periodic curve. By Lemma~\ref{lem:graph}, $Y$ is the graph of a polynomial automorphism $\sigma$ satisfying
$\sigma\circ f^r=f^{-r}\circ\sigma$. Moreover, $\Delta_{C_{\bar k}}\subset Y$ implies $C_{\bar k}\subset\Fix(\sigma)_{\bar k}$.

Observe that any polynomial automorphism commuting with $f^r$
and fixing $C_{\bar k}$ pointwise must be the identity: otherwise $f$ would admit a periodic curve.
Applying this observation to $\sigma^2$ gives $\sigma^2=\id$.
Moreover, $\sigma\neq\id$, since otherwise $f^{2r}=\id$.

For every $\tau\in\gal(\bar k/k)$, the conjugate $\sigma^\tau$ also fixes $C_{\bar k}$ pointwise and satisfies
$\sigma^\tau\circ f^r=f^{-r}\circ\sigma^\tau$.
Hence $\sigma^\tau\circ\sigma^{-1}$ commutes with $f^r$ and fixes $C_{\bar k}$ pointwise.
Thus $\sigma^\tau=\sigma$ by the above observation, which implies that $\sigma$ is indeed defined over $k$.

By \cite[Propositions~4--6]{GM03}, a nontrivial involution whose fixed locus contains a curve is conjugate over $\bar k$ to $(x,y)\mapsto(y,x)$. In particular, its fixed locus is geometrically irreducible. Hence $\Fix(\sigma)=C$.
\end{proof}

\section{Proportionality of Green functions over a number field}
\subsection{Dynamical preliminaries}
\subsubsection{Green functions and B\"ottcher coordinates}
Let $k$ be a field of characteristic zero. A \emph{generalized H\'enon map} $f:\Abb^2_k\to\Abb^2_k$ defined over $k$ is a finite composition of maps
\[
    H_i(x,y)=(y,p_i(y)-\delta_i x),
    \qquad \deg p_i=d_i\geq2,\quad \delta_i\in k^\times.
\]
Its degree is $d=\prod_i d_i$ and its Jacobian is $\prod_i\delta_i$.
The rational extensions of $f$ and $f^{-1}$ to $\Pbb^2_k$ have indeterminacy points
\[
    I_+=[1:0:0],\qquad I_-=[0:1:0],
\]
respectively. The point $I_-$ is an attracting fixed point of $f$, and $I_+$ is an attracting fixed point of $f^{-1}$. By Friedland--Milnor~\cite{FM89}, every polynomial automorphism of H\'enon type is conjugate to a generalized H\'enon map over an algebraic closure.

\medskip

Now assume $(k,|\cdot|)$ is an algebraically closed non-trivially metrized field of characteristic zero. Write $f^{\pm}(x,y)$ as
\[f^{\pm1}(x,y)=(f_{1,\pm}(x,y),f_{2,\pm}(x,y)),\]
then it satisfies
\begin{equation}\label{eq:standard-form}
\begin{aligned}
    d_1^+&:=\deg f_{1,+}=\deg_y f_{1,+}<d,
    & f_{2,+}(x,y)&=c_+y^d+q_+(x,y),\\
    d_2^-&:=\deg f_{2,-}=\deg_x f_{2,-}<d,
    & f_{1,-}(x,y)&=c_-x^d+q_-(x,y),
\end{aligned}
\end{equation}
where $c_\pm\neq0$ and $\deg q_\pm<d$.
It follows from Bedford--Smillie \cite{BS91} in the archimedean case and from Ingram \cite{Ingram14} and Kawaguchi \cite{Kawaguchi13} in the non-archimedean case that the forward and backward Green functions
\[
    G_f^\pm:=\lim_{n\to\infty}d^{-n}\log^+\|f^{\pm n}\|
\]
exist and are continuous and nonnegative, and satisfy $G_f^\pm\circ f^{\pm1}=dG_f^\pm$.
Here $\|(x,y)\|=\max\{|x|,|y|\}$. For $R$ sufficiently large, put
\[
    V_R^+:=\{|y|>\max\{|x|,R\}\},
    \qquad V_R^-:=\{|x|>\max\{|y|,R\}\}.
\]
Moreover, choosing $\beta_\pm$ with $\beta_\pm^{d-1}=c_\pm$, we have
\begin{align}\label{eq-Green-at-infinity}
    G_f^+=\log|\beta_+y|+o(1)\ \text{on }V_R^+,
    \qquad G_f^-=\log|\beta_-x|+o(1)\ \text{on }V_R^-,
\end{align}
as $\|(x,y)\|\to\infty$.

\medskip

Over $\Cbb$, the construction of B\"ottcher coordinates goes back to Hubbard--Oberste-Vorth \cite{HOV94}; see Bera--Pal--Verma \cite[Proposition~2.1 and equations~(11)--(12)]{Bera-Pal-Verma} for the following formulation.
\begin{prop}\label{prop:Basic-Green-bottcher}
For $R$ sufficiently large, there are nonvanishing holomorphic functions $\phi_f^\pm:V_R^\pm\to\Cbb$ such that
\begin{enumerate}
    \item $\phi_f^\pm\circ f^{\pm1}=(\phi_f^\pm)^d$ and $G_f^\pm=\log|\phi_f^\pm|$ on $V_R^\pm$;
    \item $\phi_f^+(x,y)\sim\beta_+y$ as $\|(x,y)\|\to\infty$ in $V_R^+$;
    \item $\phi_f^-(x,y)\sim\beta_-x$ as $\|(x,y)\|\to\infty$ in $V_R^-$.
\end{enumerate}
\end{prop}

\medskip

Finally, we close this section by recalling the following recent rigidity theorem of Cantat--Dujardin~\cite[Corollary~2.8]{CD24}.
\begin{thm}\label{thm:Cantat-Dujardin}
    Let $f$ be a complex generalized H\'enon map. Let $\Gamma$ be an embedded holomorphic disk. Suppose $\ddc(G^+_f|_\Gamma)$ non-vanishing. Then its support is not a smooth real curve.
\end{thm}

\subsubsection{Dynamical divisors}\label{sect:dynamical-divisors}

Let $f$ be a generalized H\'enon map over a field $k$ of characteristic zero, and let $C\subset\Abb^2_k$ be an algebraic curve. Let $\hat C$ be the normalization of its projective closure, and let $\iota:\hat C\to\Pbb^2$ be the induced morphism. Write $E:=\hat C\setminus C$, with $C$ understood as its inverse image in $\hat C$. For $n\geq0$, put
\[
    D_n^\pm:=(f^{\pm n}\circ\iota)^*L_\infty,
    \qquad \mathrm{where} \ \  L_\infty:=\Pbb^2\setminus\Abb^2.
\]
By \cite[Proposition~4.2]{BS91}, there exists $N_0\geq1$ such that
\[
    (f^{N_0}\circ\iota)(E)=\{I_-\},
    \qquad (f^{-N_0}\circ\iota)(E)=\{I_+\}.
\]
\begin{lem}\label{lem:stabilise-D}
For every $n\geq0$, we have $D_{N_0+n}^\pm=d^nD_{N_0}^\pm$.
\end{lem}
\begin{proof}
Put $U:=\Pbb^2\setminus\{I_+\}$. The morphism $f:U\to U$ satisfies $(f|_U)^*(L_\infty|_U)=dL_\infty|_U$. Since $(f^{N_0}\circ\iota)(\hat C)\subset\Abb^2\cup\{I_-\}\subset U$, we obtain
\[
    D_{N_0+n}^+
    =(f^{N_0}\circ\iota)^*(f|_U^n)^*(L_\infty|_U)
    =d^nD_{N_0}^+.
\]
The proof for $D_{N_0+n}^-$ is identical, using $f^{-1}$.
\end{proof}

\subsection{Arithmetic preliminaries}

\subsubsection{Adelic line bundles on curves, equidistribution and arithmeric hodge index theorem.} \label{sec:adelic}
Let $\Kbb$ be a number field, $M_\Kbb$ its set of places, and
$X$ a smooth projective geometrically irreducible curve over $\Kbb$. For $v\in M_\Kbb$, denote by $\Kbb_v$ its completion at $v$. Then let $\Cbb_v$ be the completion of an algebraic closure of $\Kbb_v$, which is complete and algebraic closed. Write $X_v^{\an}:=(X\times_\Kbb\Cbb_v)^{\an}$ for the Berkovich analytification at non-archimedean places and for the usual complex analytification at archimedean places. For a place $w$ of any number field $\Kbb'$, put $n_w:=[\Kbb'_w:\mathbb Q_w]/[\Kbb':\mathbb Q]$, where $\mathbb Q_w$ is the completion of $\mathbb Q$ at the place below $w$.

\medskip

An adelic metrized $\Qbb$-line bundle $\bar L=(L,\{\|\cdot\|_v\}_{v\in M_\Kbb})$ consists of a
$\Qbb$-line bundle $L$ with continuous local metrics induced by a common model at all but finitely many places. It is \emph{semipositive} if its metrics are uniform limits of smooth metrics with nonnegative curvature at archimedean places and of semipositive model metrics at non-archimedean places. A difference of two semipositive adelic line bundles is called \emph{integrable}. These notions are introduced by Zhang~\cite{Zhang-positive,Zhang95}. For a semipositive adelic line bundle $\bar L$, denote by $c_1(\bar L)_v$ its curvature measure on $X_v^{\an}$. This has total mass $\deg(L)$. For its construction at non-archimedean places, see Chambert-Loir~\cite{ACL} and Thuillier~\cite{Thuillier}.

\medskip

For $x\in X(\overline{\Kbb})$, put $\Kbb'=\Kbb(x)$, choose $m\geq1$ such that $mL$ is a line bundle, and let $s$ be a rational section of $mL$ that is regular and nonvanishing at $x$. The associated height is
\[
    h_{\bar L}(x):=-\frac{1}{m}
    \sum_{w\in M_{\Kbb'}}n_w\log\|s(x)\|_w.
\]
This is independent of the choices of $m,s$ and $\Kbb'$.

\begin{ex}[The standard metric on $\Pbb^1$]\label{ex:standard-metric}
\normalfont
Equip $L=\Ocal_{\Pbb^1}(1)$ with the standard metric at each place $v\in M_\Kbb$, i.e., for a linear form $s$, we define
\[
    \|s([x_0:x_1])\|_v
    =\frac{|s(x_0,x_1)|_v}{\max\{|x_0|_v,|x_1|_v\}}.
\]
The adelic line bundle $\bar L$ is semipositive. At archimedean places $c_1(\bar L)_v=\lambda_{\mathbb S^1}$ is the normalized Lebesgue measure on $\Sbb^1$.
\end{ex}

For semipositive adelic $\Qbb$-line bundles $\bar L_0,\bar L_1$ on $X$, choose $m\geq1$ such that $mL_0$ is a line bundle, and let $s$ be a nonzero rational section of $mL_0$. The induction formula of Chambert-Loir--Thuillier
\cite[Theorem~1.4]{CLT09} gives the normalized arithmetic intersection
\[
    (\bar L_0\cdot\bar L_1)
    :=\frac{1}{m}\left(
        (\bar L_1|\operatorname{div}(s))
        +\sum_{v\in M_\Kbb}n_v
\int_{X_v^{\an}}\log\|s\|_v^{-1}\,c_1(\bar L_1)_v \right),
\]
where, for $\operatorname{div}(s)=\sum_x N_x[x]$ a sum of closed points,
\[
    (\bar L_1|\operatorname{div}(s))
    =\sum_x N_x[\Kbb(x):\Kbb]h_{\bar L_1}(x).
\]
This pairing extends to integrable adelic $\Qbb$-line bundles by bilinearity.

\medskip

When $\bar{L}$ is an ample $\mathbb{Q}$-line bundle endowed with a semi-positive continuous adelic metric, following Zhang~\cite{Zhang-positive}, we have
\begin{align}
    \inf_{x\in X(\bar{\mathbb{K}})}h_{\bar{L}}(x)\leq\frac{\left(\bar{L}^2\right)}{\deg_L(X)}\leq 2 \sup_{Z\subsetneq X \ \text{finite}}\inf_{x\in (X\setminus Z)(\bar{\mathbb{K}})}h_{\bar{L}}(x).\label{eq:Zhang-ineq}
\end{align}
For such metrized line bundles, there is an automatic equidistribution of small points due to Yuan~\cite[Theorem~3.2]{Yuan} and Thuillier~\cite{Thuillier}:
\begin{thm}[Yuan~\cite{Yuan}, Thuillier~\cite{Thuillier}]\label{thm:equidistribution}
    Let $\bar{L}$ be an ample $\Qbb$-line bundle on $X$ endowed with a continuous adelic semi-positive metrization. Assume there is an infinite sequence of distinct points  $p_n\in  C(\overline{\Kbb})$ such that $h_{\bar{L}}(p_n)\to (\bar L^2)/2\deg_L(X)$. Then for any place $v\in M_\Kbb$,
    \[\frac{1}{\#\mathrm{Gal}(\overline{\Kbb}/\Kbb)\cdot p_n}\sum_{p\in \mathrm{Gal}(\overline{\Kbb}/\Kbb)\cdot p_n}\delta_p\xrightarrow[n\to\infty]{} \frac{c_1(\bar{L})_v}{\deg_L(X)},\]
    in the weak sense of measures on $X_v^\mathrm{an}$.
\end{thm}

A key result is the arithmetic Hodge index Theorem of Yuan--Zhang~\cite[Theorem~3.2]{Hodge-YZ}:
\begin{thm}[Yuan-Zhang~\cite{Hodge-YZ}]\label{thm:hodge}
If $\bar{L}_0$ and $\bar{L}_1$ are $\Qbb$-ample line bundles on $X$ with $\deg(L_0)=\deg(L_1)$ endowed with continuous adelic semi-positive metrizations, then
\[(\bar{L}_0-\bar{L}_1)^2\leq0\]
and in the equality case, there is $\bar M_0\in
\pi^*\widehat{\operatorname{Pic}}(\Spec \Kbb)_{\Qbb}$ such that $\bar{L}_1=\bar{L}_0+\bar{M}_0$.
\end{thm}

\subsubsection{Canonical heights and Green metrics of generalized H\'enon maps}\label{sec:canonical-heights}
Let $h$ denote the usual (logarithmic) Weil height on $\Pbb^2$.
\begin{thm}[Kawaguchi {\cite[Theorem~6.3]{Kawaguchi13}}]\label{thm:canonical-heights}
Let $f$ be a generalized H\'enon map of degree $d$ over a number field $\Kbb$. For every $p\in\Abb^2(\overline{\Kbb})$, the limits
\[
    \widehat h_f^\pm(p):=\lim_{n\to\infty}d^{-n}h(f^{\pm n}(p))
\]
exist and nonnegative. 

Moreover, for any finite extension $\Kbb'/\Kbb$ with $p\in\Abb^2(\Kbb')$, we have
\[
    \widehat h_f^\pm(p)
    =\sum_{w\in M_{\Kbb'}} n_w\,G_{f,w}^\pm(p).
\] 

Put $\widehat h_f:=\widehat h_f^++\widehat h_f^-$. Then $\widehat h_f(p)=0$ if and only if $p$ is $f$-periodic.
\end{thm}
The function is the \emph{canonical height} of $f$.

\medskip
Let $f$ be a generalized H\'enon map and recall the notations in $\S$~\ref{sect:dynamical-divisors}.  Put
\[
    D^\pm:=d^{-N_0}D_{N_0}^\pm,
    \qquad L^\pm:=\Ocal_{\hat C}(D^\pm).
\]
These are ample $\Qbb$-line bundles on $\hat C$.
Let $s_\pm$
be the canonical sections of $\Ocal_{\hat C}(d^{N_0}D^\pm)$. Define
\[
-\log\|s_\pm\|_v=d^{N_0}G_{f,v}^\pm
    \qquad\text{on }(\hat C\setminus E)_v^{\an}.
\]
By \cite[Proposition~1.16]{DF-MM}, this define semipositive adelic line bundles $\bar L^\pm$.
Their heights satisfy $h_{\bar L^\pm}=\widehat h_f^\pm$ on $(\hat C\setminus E)(\overline{\Kbb})$.

\medskip

For an automorphism $f$ of H\'enon type over $\Kbb$, choose a polynomial automorphism $\psi$ over a finite extension such that $g=\psi f\psi^{-1}$ is a generalized H\'enon map. Following \cite[\S4]{Kawaguchi06}, put $\widehat h_{f,\psi}:=\widehat h_g\circ\psi.$ If $\psi_1,\psi_2$ are two such conjugacies, there exist constants $c_1,c_2>0$ such that $c_1\widehat h_{f,\psi_1}\leq\widehat h_{f,\psi_2}
    \leq c_2\widehat h_{f,\psi_1}$.

\subsection{Proportionality of Green functions}

\begin{prop}\label{prop-equidistrib}
Let $f$ be a generalized H\'enon map and $C\subset\Abb^2$ a geometrically irreducible curve, both defined over a number field $\Kbb$. Suppose that there is an infinite sequence of distinct points $p_n\in C(\overline{\Kbb})$ with $\widehat h_f(p_n)\to0$. Then there exists $\alpha\in\Qbb_{>0}$ such that, for every place $v\in M_\Kbb$,
\[
    G^+_{f,v}=\alpha G^-_{f,v}
    \quad\text{on }C_v^{\an}.
\]
\end{prop}

\begin{rem}\label{rem:alpha}
    \normalfont The constant $\alpha$ is actually explicit: it is given by 
    \[\alpha=\frac{\deg(D_{N_0}^+)}{\deg(D_{N_0}^-)},\]
    see Lemma~\ref{lem:stabilise-D}.
\end{rem}

\begin{proof}
Let $\bar L^\pm$ be the semipositive adelic line bundles constructed in Section~\ref{sec:canonical-heights}.By the equidistribution Theorem~\ref{thm:equidistribution}, there is a constant $\alpha\in \mathbb{Q}^*_+$ such that, for any place $v\in M_{\mathbb{K}}$, 
\[c_1(\bar{L}^+)_v=\alpha\cdot c_1(\bar{L}^-)_v\]
as measures on $C_v^{\mathrm{an}}$. This implies in particular that $\deg(L^+)=\alpha\cdot \deg(L^-)$ and 
\begin{align*}
(\bar{L}^+-\alpha \bar{L}^-)^2 & =(\bar{L}^+)^2+\alpha^2(\bar{L}^-)^2-2\alpha (\bar{L}^+\cdot\bar{L}^-)=-2\alpha (\bar{L}^+\cdot\bar{L}^-),
\end{align*} 
where we used again that there is an infinite set of $f$-periodic points in $C$ to get $(\bar{L}^\pm)^2=0$. By construction, $f$-periodic points are points of small height for the height function $h_{\bar{L}^++\alpha\cdot \bar{L}^-}$, whence
\[0=(\bar{L}^++\alpha \bar{L}^-)^2=(\bar{L}^+)^2+\alpha^2(\bar{L}^-)^2+2\alpha(\bar{L}^+\cdot \bar{L}^-)=2\alpha(\bar{L}^+\cdot \bar{L}^-).\]
We thus have proved that
\begin{align*}
(\bar{L}^+-\alpha \bar{L}^-)^2 & =0.
\end{align*} 
By the equality case of the arithmetic Hodge index Theorem~\ref{thm:hodge}, this implies $\bar{L}^+=\alpha \cdot\bar{L}^-+\bar{M}_0$, where $\bar M_0\in
\pi^*\widehat{\operatorname{Pic}}(\Spec \Kbb)_{\Qbb}$. This means that there exists a rational function $\phi$ and a rational number $\beta>0$ such that for all $v\in M_\Kbb$,
\begin{align}\label{eq:Green_rational}
    G^+_{f,v}=\alpha G^-_{f,v} + \beta\log|\phi |_v .
\end{align}
It remains to show that $\phi$ is constant. Indeed, if $\phi$ is constant, for any $v\in M_\Kbb$, we have proved $\mathrm{supp}(dd^cG_{f,v}^+|_{C_v^\mathrm{an}})=\mathrm{supp}(dd^cG_{f,v}^-|_{C_v^\mathrm{an}})$, there is $p\in C_v^\mathrm{an}$ where $G_{f,v}^+(p)=G_{f,v}^-(p)=0$. Whence \eqref{eq:Green_rational} implies $\log|\phi(p)|_v=0$.

We now prove that $\phi$ is constant. Fix $v\in M_\Kbb$. Since $G_{f,v}^-\geq0$ on $\mathbb{A}^{2,\mathrm{an}}_{v}$, \eqref{eq:Green_rational} implies
\[\log\max\{|\phi|_v,1\}\leq\frac{1}{\beta}G_{f,v}^+ \quad \text{on} \ C_v^\mathrm{an}.\]
 Denote by $\bar E$ the pullback by $\phi$ of the standard Weil metrization of $\Ocal_{\Pbb^1}(1)$ (Example~\ref{ex:standard-metric}), as in \S~\ref{sec:adelic}. The above implies that 
\[h_{\bar{E}}\leq \frac{1}{\beta}\cdot \widehat{h}_f^+ \quad \text{on}\ \ C(\overline{\Kbb}).\]
\begin{comment}
    so that its curvature form is $\phi^*(\lambda_{\mathbb{S}^1})$ at an archimedean place, where $\lambda_{\mathbb{S}^1}$ is the normalized Lebesgue measure of the unit circle $\mathbb{S}^1.$
\end{comment}
 A sequence of distinct points $p_n\in C(\overline{\Kbb})$ with $\hat{h}_f(p_n)\to0$ is small for $\bar L^+$ and $\bar E$ at the same time. Indeed, if $p_n\in C(\overline{\mathbb{K}})$ is a sequence with $\widehat{h}_f(p_n)\to0$, we have
 \[h_{\bar{E}}(p_n)\leq \frac{1}{\beta}\widehat{h}_f^+(p_n)\leq \frac{1}{\beta}\widehat{h}_f(p_n),\]
 thus $\max\{\widehat{h}_f^+(p_n),h_{\bar{E}}(p_n)\}\to 0$ as $n\to\infty$. Fix now an archimedean place $v_0$.
Applying a second time the equidistribution theorem, we get that the measures $\ddc G^+_{f,v_0}|_C$ and $\phi^*(\lambda_{\mathbb{S}^1})$ are proportional. However, by Theorem~\ref{thm:Cantat-Dujardin}, the support of $\ddc G^+_{f,v_0}|_C$ can not be locally a smooth arc, whence a contradiction. Therefore $\phi$ is a constant. This concludes the proof. 
\end{proof}

\section{An analytic relation on B\"ottcher coordinates}

\subsection{B\"ottcher relation over a number field}\label{sect:bottcher}

\begin{prop}\label{prop-bottcher}
Let $f$ be a generalized H\'enon map and $C\subset\Abb^2_{\Kbb}$ a geometrically irreducible curve, both defined over a number field $\Kbb$. Suppose that there is an infinite sequence of points $p_n\in C(\overline{\Kbb})$ with $\widehat h_f(p_n)\to0$.
Then there exist integers $A,B\geq1$ such that, for every $t_0\in\hat C\setminus C$,
\[
    \left(\phi_f^+\circ f^{N_0}\right)^A
    =\left(\phi_f^-\circ f^{-N_0}\right)^B
\]
on a punctured neighborhood of $t_0$ in $\hat C(\Cbb)$.
\end{prop}
\begin{proof}   
    Put $E:=\hat C\setminus C$. Up to taking a finite extension of $\Kbb$, we may assume that $\beta_\pm\in\Kbb^\times$ and that every point of $E$ is $\Kbb$-rational. Write $\alpha=r/s$, where $r,s$ are coprime positive integers and $\alpha$ is the constant in Proposition~\ref{prop-equidistrib}.
The first point of Proposition~\ref{prop:Basic-Green-bottcher}, implies that the function 
\[
\frac{\left(\phi_f^+\circ f^{N_0}\right)^s}{\left(\phi_f^-\circ f^{-N_0}\right)^r}
\]
defined on a pointed neighborhood of $t_0$ in $\hat{C}(\mathbb{C})$ has constant modulus 1, hence it is a constant $\zeta_{t_0}\in \mathbb{S}^1$. We claim that $\zeta_{t_0}$ is a root of unity.

To prove it, we define the rational function 
\[
\varPhi:=\frac{\left(\beta_+\cdot \pi_2\circ f^{N_0}\right)^s}{\left( \beta_-\cdot \pi_1\circ f^{-N_0}\right)^r} \in \Kbb(\hat{C})^\times,
\]
where $\pi_i:\mathbb{A}^2\to\mathbb{A}^1$ is the projection onto the $i$ variable. By points (2) and (3) of Proposition~\ref{prop:Basic-Green-bottcher}, as $t\in C(\mathbb{C})$ converges to $t_0$, we have $\varPhi(t)\to \zeta_{t_0}$. Therefore the constant $\zeta_{t_0}=\varPhi(t_0)\in \Kbb^\times$ is algebraic.
Pick any place $v\in M_\mathbb{K}$ of $\Kbb$. By definition of $\varPhi$, we have
\begin{align*}
    \log |\varPhi(t)|_v=s\log|\beta_+\cdot\pi_2\circ f^{N_0}(t)|_v - r\log |\beta_-\cdot\pi_1\circ f^{-N_0}(t)|_v.
\end{align*}
Then the Green function asymptotics in Eq.~\eqref{eq-Green-at-infinity} and Proposition~\ref{prop-equidistrib} give
\[
\log |\varPhi(t)|_v= sd^{N_0} G^+_{f,v}(t) - rd^{N_0} G^-_{f,v}(t) +o(1)=o(1).
\]
Thus $|\varPhi(t_0)|_v=1$ for every place $v\in M_\mathbb{K}$. By Kronecker's Theorem, $\zeta_{t_0}$ is a root of unity. 

Since $E$ is finite, there exists $m\geq1$ such that $\zeta_{t_0}^m=1$ for every $t_0\in E$. Put $A:=ms$ and $B:=mr$ shows the conclusion.
\end{proof}

\subsection{Preparation for specialization: spreading out}\label{sect:spreading-out}
Let $f$ be a generalized H\'enon map and let $C\subset\mathbb{A}^2_\mathbb{C}$ be an irreducible curve that contains infinitely many periodic points of $f$. Choose a field $K\subset\Cbb$, finitely generated over $\bar{\Qbb}$, over which $f$ and $C$ are defined.

Let $\hat C_K$ be the normalization of the projective closure of $C_K$ and put $E_K:=\hat C_K\setminus C_K$. Up to taking a finite extension of $K$ inside $\Cbb$, we may assume that all points of $E_K$ are $K$-rational and that the constants $\beta_\pm$ of Proposition~\ref{prop:Basic-Green-bottcher} belong to $K^\times$.

Let $R\subset K$ be a finitely generated $\bar{\mathbb{Q}}$-algebra whose fraction field is $K$ and put $S=\Spec (R)$. For any $s\in S$, $f_s$ and $C_s$ will be denoted as the map and curve obtained by specialization at $s$.
By Proposition~\ref{prop-bottcher}, we can reduce to the case when $\dim S>0$. The embedding $K\hookrightarrow\Cbb$ defines a point $s_0\in S(\Cbb)$ with $f_{s_0}=f$ and $C_{s_0}=C$ whose image in $S$ is the generic point.

After shrinking $S$, we may assume that  $S$ is smooth and $C_s$ is geometrically irreducible. There exists a smooth projective family $\hat{\mathcal C}\to S$ whose fibers are the normalizations of the projective closures of the curves $C_s$. Write $\CCal\subset \hat{\CCal}$ for the open subset corresponding to $C_s$. 

After shrinking $S$ again, the points of $E_K$ extend to disjoint sections $\xi_1,\ldots,\xi_\ell$ with
\begin{align*}
    \Ecal:=\hat \CCal \setminus \CCal=\bigsqcup_{i=1}^{\ell}\xi_i(S).
\end{align*}
Moreover, we have
\[
    f_s^{N_0}(\xi_i(s))=I_-,
    \qquad f_s^{-N_0}(\xi_i(s))=I_+
\]
for all $s\in S$ and all $i$, and 
\[
D_{N_0,s}^\pm=\sum_{i=1}^{\ell}m_i^\pm\xi_i(s),
    \qquad m_i^\pm>0,
\]
where the positive integers $m_i^\pm$ are independent of $s$. 

Finally, after shrinking $S$, we may assume that $\beta_\pm\in \Ocal(S)^\times$ are regular and nowhere zero on $S$ and use their values at $s$ to choose the B\"ottcher coordinates of $f_s$.

\subsection{B\"ottcher relation over $\Cbb$}\label{sect:bottcher-complex}
\begin{prop}\label{prop:complex-bottcher}
Let $f$ be a generalized H\'enon map of $\Cbb^2$ and let $C\subset\Abb^2_\Cbb$ be an irreducible curve containing infinitely many $f$-periodic points. Then there exist a point $t_0\in\hat C\setminus C$ and integers $A,B\geq1$ such that
\[
    \left(\phi_f^+\circ f^{N_0}\right)^A
       =\left(\phi_f^-\circ f^{-N_0}\right)^B
\]
on a punctured neighborhood of $t_0$ in $\hat C(\Cbb)$.
\end{prop}

\begin{proof}
For every $s\in S(\bar{\Qbb})$, the map $f_s$ and the curve $C_s$ are defined over a number field. By \cite[Lemma~5.3]{DF-MM}, the curve $C_s$ contains infinitely many $f_s$-periodic points. Proposition~\ref{prop-equidistrib} and Remark~\ref{rem:alpha} therefore give
\[
    G_{f_s}^+=\alpha G_{f_s}^-\quad\text{on }\ C_s(\Cbb),
    \qquad
    \alpha=\frac{\deg D_{N_0,s}^+}{\deg D_{N_0,s}^-}
           =\frac{\sum_i m_i^+}{\sum_i m_i^-}.
\]
In particular, $\alpha\in\Qbb_{>0}$ does not depend on $s$.

The Green functions are continuous in both the parameter and the space variables (see \cite[Proposition 3.4]{BS91}). Thus 
\begin{equation}\label{eq:complex-green-family}
    G_{f_s}^+=\alpha G_{f_s}^-\quad\text{on }\ C_s(\Cbb)
\end{equation}
for every $s\in S(\Cbb)$.

Fix any section, say $\xi:=\xi_1$ and put $a:=\sum_i m_i^-$ and $b:=\sum_i m_i^+$, so that $\alpha=b/a$. For every $s\in S(\Cbb)$, the function
\[
q_s:=\frac{\left(\phi_{f_s}^+\circ f_s^{N_0}\right)^a}{\left(\phi_{f_s}^-\circ f_s^{-N_0}\right)^b}
\]
is holomorphic on a small punctured neighborhood of $\xi(s)$ in $\hat C_s$ and has modulus $1$. Hence it is a constant $\lambda_s\in\mathbb S^1$.

Let $\pi_1,\pi_2:\Abb^2\to\Abb^1$ be the coordinate projections and consider the rational function
\[
\varPhi:=\frac{\left(\beta_+\cdot\pi_2\circ f^{N_0}\right)^a}{\left(\beta_-\cdot\pi_1\circ f^{-N_0}\right)^b}
\]
on $\hat{\mathcal C}$. After shrinking $S$, we may assume that $\varPhi$ is regular and nonzero on a neighborhood of $\xi(S)$. Thus $\lambda:=\xi^*\varPhi\in\Ocal(S)^\times$. We infer that
\[
\lambda(s)=\varPhi_s(\xi(s))=\lambda_s
\]
for every $s\in S(\Cbb)$. Therefore $\lambda$ is not an open map which implies that $\lambda$ is constant.

Choose $s_1\in S(\bar{\Qbb})$. The proof of Proposition~\ref{prop-bottcher}, applied to $(f_{s_1},C_{s_1})$ at $\xi(s_1)$, shows that $\lambda(s_1)=\lambda$ is a root of unity. Let $m\geq1$ be its order. Put $t_0:=\xi(s_0)$, $A:=ma$ and $B:=mb$. Then $q_{s_0}^m=\lambda^m=1$ gives,
\begin{equation}\label{eq:complex-boundary-phase}
    \left(\phi_f^+\circ f^{N_0}\right)^A
       =\left(\phi_f^-\circ f^{-N_0}\right)^B,
\end{equation}
on a punctured neighborhood of $t_0$ in $\hat C(\Cbb)$.
\end{proof}

\section{A relation on B\"ottcher coordinates implies non-Zariski density}\label{sect:xie}

In this section, we prove the following and deduce the main results of the paper.

\begin{thm}\label{tm-non-density-from-Bottcher}
Let $f$ be a generalized H\'enon map of $\Cbb^2$ and let $C\subset\Abb^2_\Cbb$ be an irreducible curve. Assume that there exist a point $t_0\in\hat C\setminus C$ and integers $A,B\geq1$ such that
\begin{align}\label{eq:AB}
\left(\phi_f^+\circ f^{N_0}\right)^A
       =\left(\phi_f^-\circ f^{-N_0}\right)^B    
\end{align}
on a punctured neighborhood of $t_0$ in $\hat C$. Then $\bigcup_{n\in\Zbb}F^n(\Delta_C)$ is not Zariski dense in $\Cbb^4$.
\end{thm}

\subsection{Xie's algebraicity criterion} 
Let $K$ be a finite extension of $k(t)$, where $k$ is a  field of characteristic $0$ and let $M_K$ be the set of places of the field $K$. Let $X$ be a smooth projective variety over $K$. Recall the following result of Xie (see~\cite[Theorem~1.4]{Xie2024}). 
\begin{thm}[Xie]\label{tm-Xie}
 Let $h : X (K ) \to\mathbb{R}$ be any Weil height. Let $v \in M_K$ and let $\mathcal{V}$ be a
local analytic subspace of $X_v$ at a point $o$ which is defined over $K$. Let $(x_n)_{n\geq0}$ be a sequence of $K$-points in $X(K)$ such that:
\begin{enumerate}
    \item $x_n\in \mathcal{V}$ for all $n\geq0$,
    \item $(x_n)_n$ converges to $o$ in $X(K_v)$ as $n\to\infty$,
    \item $\liminf_{n\to\infty}\frac{-\log d_v(x_n,o)}{h(x_n)}>0$.
\end{enumerate}
Set $H := \bigcap_{m\geq0}\overline{\{x_n ,n \geq m\}}^{\mathrm{Zar}} \subseteq X$ . Then $\dim H \le \dim (\Vcal)$.
\end{thm}
Here $d_v$ stand for the distance inherited from the chordal distance on $\mathbb{P}^N(K_v)$ through an embedding $X\hookrightarrow\mathbb{P}^N$ and the chordal distance is defined as
\[d_v(x,y):=\frac{\max_{i,j}|X_iY_j-X_jY_i|}{\|X\|\cdot\|Y\|},\]
where $x=[X_0:\cdots:X_N]$ and $y=[Y_0:\cdots:Y_N]$ as points of $\mathbb{P}^N(K_v)$. Note that, in the case when $h(x_n)\to\infty$ as $n\to\infty$, condition (2) is a consequence of condition (3).

 In particular, if we take $X=\Pbb^2\times \Pbb^2$, then using the Segre embedding $X\hookrightarrow \Pbb^{8}$, we get the following formula. Take $o:=[I_-:I_+]$ and $X:=[X_0:1:X_2]$ and $Y:=[1:Y_1:Y_2]$ such that $|X_0|_v,|X_2|_v,|Y_1|_v,|Y_2|_v\le 1$. Then 
\begin{align*}
    d_v((X,Y),o)&=\max \{|X_0|_v,|X_2|_v,|Y_1|_v,|Y_2|_v \}\\
    &=\max\left\{
\max\{|X_0|_v,|X_2|_v\},
\max\{|Y_1|_v,|Y_2|_v\}
\right\}.
\end{align*}
We infer that
\begin{align}\label{eq:distance}
    d_v((X,Y),o)=\max\{d_{\mathbb{P}^2,v}(X,I_-),d_{\mathbb{P}^2,v}(Y,I_+)\}.
\end{align}

\subsection{Proof of Theorem~\ref{tm-non-density-from-Bottcher}}
Let $K$ be the function field $\mathbb{C}(\hat{C})$ and $v\in M_K$ be the place of $K$ induced by the norm $|\cdot|_t$. Let $X:=\mathbb{P}^2\times\mathbb{P}^2$ and $o=(I^-,I^+)$. 

We first construct the analytic germ used in the proof. By the removable singularity theorem, the functions $\psi_f^\pm:=1/\phi_f^\pm$ extend holomorphically to neighborhoods of $I_\mp$, respectively, and satisfy $\psi_f^+(I_-)=\psi_f^-(I_+)=0$. Indeed, near $I_-$ the coordinates $(u,z)=(x/y,1/y)$ give $\psi_f^+\sim z/\beta_+$; the argument near $I_+$ is the same. Thus
\[
    \mathcal V_f:=\{(p,q):\psi_f^+(p)^A=\psi_f^-(q)^B\}
\]
defines a proper complex analytic hypersurface germ at $o$. Its defining Taylor series has coefficients in $\Cbb$, on which $v$ is trivial, and therefore also defines a $v$-adic analytic hypersurface $\mathcal V$ near $o$ in $X(K_v)$.

Let $c\in C(K)$ be the generic point of $C$. Denote again by $F$ the rational extension to $X_K$ of $(f,f^{-1})$, and put
\[
    x_n:=F^{N_0+n}(c,c)\in X(K)
\]
for every $n\geq0$.
Since 
\[
    \psi_f^+\bigl(f^{N_0+n}(c)\bigr)^A
       =\psi_f^-\bigl(f^{-(N_0+n)}(c)\bigr)^B,
\]
we have that $x_n\in \mathcal{V}$. 

Let $h:X(K)\to\mathbb{R}_+$ be a Weil height associated to the line bundle $\Lcal:=\Ocal(1,1)$. 
\begin{lem}\label{lem-height}
    There is a constant $C>0$ such that $h(x_n)=C\cdot d^n+O(1)$.
\end{lem}
\begin{proof}
    For every $n\geq0$, put
\[
    \gamma_n:=\bigl(f^{N_0+n}\circ\iota,f^{-(N_0+n)}\circ\iota\bigr):\hat C\to X.
\]
Then,
\[
    h(x_n)=\deg\gamma_n^*\mathcal L+O(1)
           =\deg D_{N_0+n}^++\deg D_{N_0+n}^-+O(1).
\]
By Lemma~\ref{lem:stabilise-D}, this equals
\[
    d^n\bigl(\deg D_{N_0}^++\deg D_{N_0}^-\bigr)+O(1).
\]
Taking $C:=\deg D_{N_0}^++\deg D_{N_0}^->0$ proves the lemma.
\end{proof}
\begin{lem}\label{lem-distance}
There is a constant $C'>0$ such that, for all $n\geq0$,
\[-\log d_v(x_n,o)= C'\cdot d^n.\]
\end{lem}
\begin{proof}
Let $t$ be a local coordinate on $\hat{C}$ around $t_0$ fixed above. Recall that for any rational map $g\in \Cbb(\hat{C})$, we have $|g|_v=\exp(-\mathrm{ord}_v(g))$. Write 
\[ 
(p_n^+,q_n^+)=f^{N_0+n}\circ\iota \qquad \mathrm{and}\qquad (p_n^-,q_n^-)=f^{-(N_0+n)}\circ\iota.
\]
Then $x_n=([p_n^+:q_n^+:1],[p_n^-:q_n^-:1])$ and locally at $t_0$, we can write
\[
p_n^+(t)=t^{a_n^+}u_n(t),\qquad
q_n^+(t)=t^{b_n^+}v_n(t),
\]
where $u_n(0),v_n(0)\ne0$. By the choice of $N_0$, we have $b_0^+:=\operatorname{ord}_{t_0}(q_0^+)<0$. Since $(f^{N_0+n}\circ\iota)(t_0)=I_-$, the coefficient of $t_0$ in $D_{N_0+n}^+$ is $-b_n^+$. Lemma~\ref{lem:stabilise-D} therefore gives $b_n^+=d^n b_0^+$.

Moreover, the leading $y^{d_1^+}$-term of $f_{1,+}$ dominates near $I_-$. Thus \eqref{eq:standard-form} gives $a_n^+=d_1^+d^{n-1}b_0^+$ for $n\geq1$, and
\begin{align*}
d_v([p_n^+:q_n^+:1],I_-)&=d_v([p_n^+/q^+_n:1:1/q^+_n],I_-)\\
    &=\max\{\exp(b^+_n-a^+_n),\exp(b^+_n) \}\\
    &=\exp(b^+_n-a^+_n)=\exp((d-d_1^+)d^{n-1}b_0^+).
\end{align*}

Similarly, if
\[
p_n^-(t)=t^{a_n^-}w_n(t),\qquad
q_n^-(t)=t^{b_n^-}z_n(t),
\]
where $w_n(0),z_n(0)\ne0$, we get
\begin{align*}
     d_v([p_n^-:q_n^-:1],I_+)=\exp((d-d_2^-)d^{n-1} a_0^{-} ).
\end{align*}

By \eqref{eq:distance}, we have proved that
\begin{align*}
    -\log d_v(x_n,o)&=-\log \max \{d_v([p_n^+:q_n^+:1],I_-),  d_v([p_n^-:q_n^-:1],I_+)\}\\
&=d^n \min\{(d_1^+-d)b_0^+,(d_2^--d)a_0^- \}. 
\end{align*}
We concludes the proof letting $C':=\min\{(d_1^+-d)b_0^+,(d_2^--d)a_0^- \}>0$.
\end{proof}

The combination of Lemma~\ref{lem-height} and Lemma~\ref{lem-distance} gives
\[\lim_{n\to\infty}\frac{-\log d(x_n,o)}{h(x_n)}= \lim_{n\to\infty}\frac{C'}{C+O(d^{-n})}=\frac{C'}{C}>0.\]
Moreover, Lemma~\ref{lem-height} also guaranties that $h(x_n)\to\infty$ as $n\to\infty$, proving that assumption (3) (and thus assumption (2)) of Theorem~\ref{tm-Xie} are satisfied by the sequence $(x_n)_n$. Thus, if we put $H:=\bigcap_{j\ge1}\overline{\{x_n:n\ge j\}}^{\mathrm{Zar}}$, where the Zariski closure is taken in $\Abb^4$, then $\dim H\le 3$. Moreover, by construction, $F(H)=H$. Applying Lemma~\ref{lem:dim_2}, we infer that $\dim H =2$.

Let $\Hcal$ be the Zariski closure of $H$ in $\hat{C}\times \Abb^4$. It has dimension $\dim \Hcal=\dim H +\dim \hat{C}=3$. Denote by $\mathrm{pr}_2$ the second projection $\hat{C}\times \Abb^4_\Cbb\to \Abb^4_\Cbb$. Then $Y:=\mathrm{pr}_2(\Hcal)$ is $F$-fixed and has dimension at most three. By construction of $H$, $F^n(\Delta_C)\subset Y$ for all sufficiently large $n$. Applying $F^{-1}$, we get $\Delta_C\subset Y$. It follows then $\overline{\bigcup_{n\in\Zbb}F^n(\Delta_C)}^{\mathrm{Zar}} \subsetneq\Abb^4_{\Cbb}$.

\subsection{Proofs of the main results}
We are now in position to prove Theorem~\ref{thm:main}.

\begin{proof}[Proof of Theorem~\ref{thm:main}]
    For the direct implication, since the number of periodic points on a curve does not change by conjugacy, we may assume that $f$ is a generalized H\'enon map. Apply Proposition~\ref{prop:complex-bottcher}, Theorem~\ref{tm-non-density-from-Bottcher} and then Lemma~\ref{lem:keylemma}. The inverse implication is due to Dujardin--Favre~\cite[Proposition 7.1]{DF-MM}.
\end{proof}

Over number fields, we have the stronger following Bogomolov-type version:

\begin{thm}\label{mainthm:nbfield}
       Let $f$ be a polynomial automorphism of $\Abb^2$ of H\'enon type. Let $C\subset \Abb^2$ be a geometrically irreducible algebraic curve, both defined over a number field. 
    Then there is an infinite sequence of distinct points $p_n\in C(\overline{\Kbb})$ with $\widehat{h}_f(p_n)\to0$  if and only if there exists an integer $r\ge 1$ and a polynomial involution $\sigma$ of $\Abb^2$ defined over $\Kbb$ such that its fixed points are $C$ and $\sigma \circ f^{r} \circ \sigma ^{-1}=f^{-r}$.
\end{thm}

\begin{proof}
For the direct implication, choose a polynomial automorphism $\psi$ over a finite extension of $\Kbb$ such that $g:=\psi f\psi^{-1}$ is a generalized H\'enon map. By Section~\ref{sec:canonical-heights},
$\widehat h_g(\psi(p_n))\to0$. Proposition~\ref{prop-bottcher} and
Theorem~\ref{tm-non-density-from-Bottcher}, applied to $(g,\psi(C))$, show that
the orbit of $\Delta_{\psi(C)}$ under $(g,g^{-1})$ is not Zariski dense. Applying $\psi^{-1}\times\psi^{-1}$ gives the non Zariski density for the orbit of $\Delta_C$ under $(f,f^{-1})$. Lemma~\ref{lem:keylemma}, applied to $(f,C)$ over
$\Kbb$ gives the required involution defined over $\Kbb$. The inverse implication is again due to Dujardin--Favre~\cite[Proposition 7.1]{DF-MM}.
\end{proof}

\bibliographystyle{alpha}
\bibliography{Mybio}
\end{document}